\documentclass[11pt]{article}

\usepackage{amsfonts}

\usepackage{amsmath}

\newcommand{\beq}{\begin{equation}}

\newcommand{\eeq}{\end{equation}}

\newcommand{\bea}{\begin{eqnarray}}

\newcommand{\eea}{\end{eqnarray}}

\newcommand{\beas}{\begin{eqnarray*}}

\newcommand{\eeas}{\end{eqnarray*}}

\newtheorem{theorem}{Theorem}[section]

\newtheorem{proposition}[theorem]{Proposition}

\newtheorem{corollary}[theorem]{Corollary}

\newtheorem{remark}[theorem]{Remark}

\newtheorem{example}[theorem]{Example}

\newtheorem{examples}[theorem]{Examples}

\newtheorem{foo}[theorem]{Remarks}

\newenvironment{proof}{\addvspace{\medskipamount}\par\noindent{\it

Proof}.}{\unskip\nobreak\hfill$\Box$\par\addvspace{\medskipamount}}

\title{Chen series and Atiyah-Singer theorem\thanks{This version of
\cite{Bau-Chen} was proofread and polished by the models Sol 6.2 and
Fable 5.2 under~the~supervision~of~the~author}}

\author{Fabrice Baudoin
\\
\\\small Department of Mathematics, Aarhus University
\\\small fbaudoin@math.au.dk}

\date{}

\begin{document}

\maketitle

\begin{abstract}
The purpose of this work is to give a new and short proof of the
Atiyah-Singer local index theorem for the Dirac operator on the
spin bundle. This proof is obtained by using heat semigroups
approximations
 based on the truncation of Brownian Chen series.

\end{abstract}

\tableofcontents

\section{Introduction}

\

The goal of this paper is to give a new and short proof of the
local Atiyah-Singer index theorem by using approximations of heat
semigroups. The heat equation approach to index theorems is not
new: It was suggested by Atiyah-Bott \cite{At-Bo} and
McKean-Singer \cite{McK-Sin}, and first carried out by Patodi
\cite{Pat} and Gilkey \cite{Gi} (see also Getzler \cite{Get}).
Bismut in \cite{Bi2} introduces
stochastic methods based on Feynman-Kac formula. For probabilistic
approaches that are mainly based on Bismut's ideas, we also refer
to \cite{Leandre} and Chapter 7 of \cite{Hsu}. For a complete
survey on (non probabilistic) heat equation methods for index
theorems, we refer to the book \cite{Ber-Ge-Ve}.

\

However, in our approach we will see that the $A$-genus appears in
a natural way,  from purely local computations on approximations
of  heat semigroups. Our method relies on explicit approximations
of the holonomy over the heat equation on vector bundles and
unlike the other probabilistic approaches does not involve the
Feynman-Kac formula nor the technics of stochastic differential
geometry.

\

The idea is the following. Let $\mathbf{P}_t$ denote a heat
semigroup. In recent works, see Baudoin \cite{Bau} and
Lyons-Victoir \cite{Ly-Vi}, by using Brownian Chen series, it has
been pointed out that $\mathbf{P}_t$ admits a formal
representation as the expectation of the exponential of a random
Lie series. The truncature of this Lie series leads to explicit
approximations of $\mathbf{P}_t$. More precisely, one gets a
family of operators $\mathbf{P}^N_t$, $N \ge 1$, such that (in the
$\sup$ norm)
\begin{align}\label{pr}
\mathbf{P}_t=\mathbf{P}_t^N + O(t^{\frac{N+1}{2}}), \quad t
\rightarrow 0.
\end{align}
This point of view has  been used in Lyons-Victoir \cite{Ly-Vi} to
provide cubature formulae on Wiener space that give efficient
numerical approximations of solutions of heat equations.

\

Assume now that $\mathbf{P}_t$ is the heat semigroup associated
with the Dirac operator on the Clifford module over a compact
$d$-dimensional spin manifold, $d$ even, and let $p_t$ denote its
kernel. The estimate (\ref{pr}) is a weak estimate: It holds on a
fixed smooth section. We will see that it can be turned into a
pointwise estimate for the kernels on the diagonal, at the price
of one more order in the truncation, and deduce
\[
\mathbf{Str} \text{ } p_t (x,x)=\mathbf{Str} \text{ } p_t^{d+1}
(x,x) + O(t^{\frac{1}{2}}), \quad t \rightarrow 0,
\]
where $\mathbf{Str}$ denotes the supertrace and $p_t^{d+1}$ the
kernel of $\mathbf{P}_t^{d+1}$. The Lie structure that explicitly
appears in $\mathbf{P}_t^{d+1}$ now leads to algebraic
cancellations that imply
\[
\lim_{t \rightarrow 0} \mathbf{Str} \text{ } p^{d+1}_t
(x,x)=\lim_{t \rightarrow 0} \mathbf{Str} \text{ } p_t^2 (x,x).
\]
Since, from McKean-Singer theorem, the integral over the manifold
of the supertrace of $p_t (x,x)$ has to be constant and equal to
the index of the Dirac operator, the local index theorem follows
from the easy computation of this last limit.

\

The paper is organized as follows. In the first part,  we survey
results on random Chen series that are needed for the construction
of approximations of heat semigroups. In the second part, we use
these series to construct explicit approximations of the holonomy
above general heat equations on vector bundles.  Finally, in the
third part, we develop the idea hinted above to provide a new
short proof of the Atiyah-Singer local index theorem for the Dirac
operator on the spin bundle.

\

\textbf{Acknowledgment:} \textit{The author would like to thank M.
Hairer for stimulating discussions}.

\section{Chen series}

We introduce here notations that will be used throughout the paper
and survey some results on Chen series that will be later needed.
Basic background on Chen series with respect to regular paths can
be found in \cite{Che} and background on Chen series with respect
to Brownian paths can be found in Chapter 1 of \cite{Bau} (see
also \cite{Fli} and \cite{Ly-Vi}). Let us note that the Chen
series with respect to Brownian paths is also called, in the rough
paths theory of Lyons, the signature of the Brownian motion.

Let $\mathbb{R} [[X_0,...,X_d]]$ be the non commutative algebra
over $\mathbb{R}$ of the formal series with $d+1$ indeterminates,
that is the set of series
\[
Y=\sum_{k \geq 0} \sum_{i_1,...,i_k} a_{i_1,...,i_k}
X_{i_1}...X_{i_k}.
\]

The exponential of $Y \in \mathbb{R} [[ X_0 ,..., X_d ]]$ is
defined by
\[
\exp (Y)=\sum_{k=0}^{+\infty} \frac{Y^k}{k!}.
\]
We define the bracket between two elements $U$ and $V$ of
$\mathbb{R} [[ X_0 ,..., X_d ]]$ by
\[
[U,V]=UV-VU.
\]
If $I=(i_1,...,i_k) \in \{ 0,..., d \}^k$ is a word, we denote by
$X_I$ the commutator defined by
\[
X_I = [X_{i_1},[X_{i_2},...,[X_{i_{k-1}}, X_{i_{k}}]...].
\]
We will denote by $\mathcal{S}_k$ the set of the permutations of
$\{1,...,k\}$.  If $\sigma \in \mathcal{S}_k$, we denote
$e(\sigma)$ the cardinality of the set
\[
\{ j \in \{1,...,k-1 \} , \sigma (j) > \sigma(j+1) \},
\]
and $\sigma (I)$ the word $(i_{\sigma(1)},...,i_{\sigma(k)})$.

Let us now consider a $d$-dimensional standard Brownian motion
$(B_t)_{t \geq 0}=(B^1_t,...,B^d_t)_{t \geq 0}$. We use the
convention that $B^0_t=t$. If $I=(i_1,...i_k) \in \{0,...,d\}^k$
is a word with length $k$, the iterated Stratonovich integral
\begin{equation*}
\int_{\Delta^k [0,t]} \circ dB^I= \int_{0 \leq t_1 \leq ... \leq
t_k \leq t} \circ dB^{i_1}_{t_1} \circ ... \circ dB^{i_k}_{t_k},
\end{equation*}
can be defined as the limit in $p$-variation, $p>2$,
\[
\lim_{ n \rightarrow +\infty} \int_{\Delta^k [0,t]}  (dB^n)^I,
\]
where $B^n$ denotes the piecewise linear interpolation of the
paths of  $(B_u)_{0 \le u \le t}$ along the dyadic subdivision of
$[0,t]$.

With the notation,
\[
\Lambda_I (B)_t= \sum_{\sigma \in \mathcal{S}_k} \frac{\left(
-1\right) ^{e(\sigma )}}{k^{2}\left(
\begin{array}{l}
k-1 \\
e(\sigma )
\end{array}
\right) } \int_{\Delta^k [0,t]} \circ dB^{\sigma^{-1}(I)},
\]
we have the following theorem:

\begin{theorem}\label{Chen-Strichartz brownien}
\[
1 + \sum_{k=1}^{+\infty} \sum_{I \in \{0,1,...,d\}^k} \left(
\int_{\Delta^k [0,t]} \circ  dB^I \right) X_{i_1} ... X_{i_k}
=\exp \left( \sum_{k \geq 1} \sum_{I \in \{0,1,...,d\}^k}\Lambda_I
(B)_t X_I \right), \text{ }t \geq 0.
\]
\end{theorem}
This theorem is due to Chen \cite{Che} and Strichartz \cite{Stri}
that prove that the above result holds for absolutely continuous
paths. The result for Brownian paths is pointed out in Fliess
\cite{Fli}. And finally, Lyons \cite{Ly}, with rough paths theory,
shows that it actually can be extended to very general paths.

If
\[
Y=\sum_{k \geq 0} \sum_{I=(i_1,...,i_k)} a_{i_1,...,i_k}
X_{i_1}...X_{i_k}
\]
is a  random series, that is if the coefficients are real random
variables defined on a probability space, we will denote
\[
\mathbb{E}(Y)=\sum_{k \geq 0} \sum_{I=(i_1,...,i_k)}
\mathbb{E}(a_{i_1,...,i_k}) X_{i_1}...X_{i_k}
\]
as soon as this expression makes sense, that is as soon as for
every $I=(i_1,...,i_k)$,
\[
\mathbb{E}(\mid a_{i_1,...,i_k} \mid)<+\infty,
\]
where $\mathbb{E}$ stands for the expectation.

\begin{theorem} \label{main theorem} (See \cite{Bau},
\cite{Ly-Vi}). We have
\[
\exp \left( t \left(X_0+\frac{1}{2}\sum_{i=1}^d X_i^2 \right)
\right)=\mathbb{E} \left( \exp \left( \sum_{k \geq 1} \sum_{I \in
\{0,1,...,d\}^k} \Lambda_I (B)_t X_I \right) \right), \quad t \ge
0.
\]
\end{theorem}

\section{Approximation of elliptic  heat kernels on vector bundles}

In the spirit of Azencott \cite{Az1}, Ben Arous \cite{Ben2},
Castell \cite{Cast} and Lyons-Victoir \cite{Ly-Vi}, we use in this
section Brownian Chen series in order to provide efficient
approximations of heat semigroups and corresponding heat kernels.
The idea is to truncate the Lie series that appear in the  formal
representation of the heat semigroup given by Theorem \ref{main
theorem}. As we shall see, this truncation that has already found
applications for cubature formulae \cite{Ly-Vi}, is also
particularly efficient to approximate the holonomy over the heat
equation in a vector bundle.

\

Let $\mathbb{M}$ be a $d$-dimensional compact smooth Riemannian
manifold and let $\mathcal{E}$ be a finite-dimensional vector
bundle over $\mathbb{M}$. We denote by $\Gamma ( \mathbb{M},
\mathcal{E} )$ the space of smooth sections. Let now $\nabla$
denote a connection on $\mathcal{E}$.

We consider the following linear partial differential equation
\begin{equation}\label{Hormander PDE}
\frac{\partial \Phi}{\partial t}=\mathcal{L} \Phi,\quad
\Phi(0,x)=f(x),
\end{equation}
where $\mathcal{L}$ is an operator on $\mathcal{E}$  that can be
written
\[
\mathcal{L}=\nabla_{0}+\frac{1}{2} \sum_{i=1}^d \nabla_{i}^2,
\]
with
\[
\nabla_i=\mathcal{F}_i +\nabla_{V_i}, \quad 0 \le i \le d,
\]
the $V_i$'s being smooth vector fields on $\mathbb{M}$ and the
$\mathcal{F}_i$'s being smooth potentials (that is sections of the
bundle $\mathbf{End}(\mathcal{E})$). It is known that the solution
of (\ref{Hormander PDE}) can be written
\[
\Phi (t,x)=(e^{t \mathcal{L}}f)(x)=\mathbf{P}_t f (x).
\]
If $I\in \{0,1,...,d\}^k$ is a word, we denote
\[
\nabla_I= [\nabla_{i_1},[\nabla_{i_2},...,[\nabla_{i_{k-1}},
\nabla_{i_{k}}]...].
\]
and
\[
d(I)=k+n(I),
\]
where $n(I)$ is the number of 0 in the word $I$.

Following this idea, we would like to consider, for $N \ge 1$, the
operators
\[
\mathbb{E} \left( \exp \left( \sum_{I, d(I) \le N } \Lambda_I
(B)_t \nabla_{I} \right) \right).
\]
Some care is however needed. First, when $N \ge 3$, the random
variables $\Lambda_I (B)_t$ can have components in Wiener chaoses
of order greater than or equal to $3$. Exponential integrability
can then fail, so that, in the presence of potentials, the above
expectation may be infinite. Secondly, when the coefficients
$\Lambda_I (B)_t$ are large, the flow of the vector field $\sum
\Lambda_I (B)_t V_I$ may come back to its starting point after a
long excursion (think of the flow of a rotation on a sphere), and
this creates singularities in the kernels that will be considered
later. Since both phenomena only concern events whose probability
is $O(t^m)$ for every $m \ge 0$, we get rid of them by a cutoff
that we now describe.

\

For $t>0$ and $b=(b_1,...,b_d) \in \mathbb{R}^d$, let us denote
\[
B^{t,b}_s=B_s +\frac{s}{t} (b-B_t), \quad 0 \le s \le t.
\]
We have $B^{t,B_t}=B$ and $B^{t,b}_s=B^{t,0}_s+\frac{s}{t} b$. The
process $B^{t,0}$ is a Brownian loop with length $t$ which is
independent of the random variable $B_t$, and $B^{t,b}$ is a
Brownian bridge from $0$ to $b$ with length $t$. The iterated
integrals of $B^{t,b}$ and the coefficients $\Lambda_I
(B^{t,b})_t$ are defined as in Section 2, with the convention that
the component with index $0$ of $B^{t,b}_s$ is $s$. Since for $1
\le i \le d$,
\[
dB^{t,b,i}_s=dB^{t,0,i}_s+\frac{b_i}{t} ds,
\]
by expanding the iterated integrals we see that for every word
$I=(i_1,...,i_k)$, $\Lambda_I (B^{t,b})_t$ is a polynomial in $b$:
\begin{align}\label{polynomial}
\Lambda_I (B^{t,b})_t=\sum_{A} \left( \prod_{a \in A}
\frac{b_{i_a}}{t} \right) \Lambda_{I^A} (B^{t,0})_t,
\end{align}
where the sum is taken over the subsets $A$ of $\{ a , i_a \neq 0
\}$ and where $I^A$ is the word obtained from $I$ by replacing by
$0$ the letters $i_a$, $a \in A$.

\

Let now $\chi: \mathbb{R}_{\ge 0} \rightarrow [0,1]$ be a smooth
function such that $\chi=1$ on $[0,1]$ and $\chi=0$ on
$[4,+\infty)$, and let $0<\rho \le 1$ be a constant. For $N \ge
1$, we consider the random variable
\[
\chi^N_t (B)=\chi \left( \frac{1}{\rho^2} \left( \mid B_t \mid^2
+\sum_{J} \frac{\Lambda_J (B^{t,0})_t^2}{t^{2n(J)}}
\right)\right),
\]
where the sum is taken over the words $J$ such that $2 \le \mid J
\mid \le N$ and that have at least one letter different from $0$.
In this definition, each letter $0$ of $J$ is compensated by a
factor $\frac{1}{t}$: By the scaling property of Brownian motion,
the family
\[
\left( B_t , \left( \frac{\Lambda_J (B^{t,0})_t}{t^{n(J)}}
\right)_J \right)
\]
has the same law as
\[
\left( \sqrt{t} B_1 , \left( t^{\frac{\mid J \mid -n(J)}{2}}
\Lambda_J (B^{1,0})_1 \right)_J \right),
\]
and $\mid J \mid -n(J) \ge 1$. Since iterated integrals have
moments of all orders, we deduce that for every $m \ge 0$,
\begin{align}\label{cutoff1}
\mathbb{P} \left( \chi^N_t (B) \neq 1 \right)=O(t^m), \quad t
\rightarrow 0.
\end{align}
On the other hand, if all the letters of a word $J$ are $0$, then
$\Lambda_J (B^{t,0})_t$ is a constant multiple of $t^{\mid J
\mid}$. From (\ref{polynomial}), we therefore see that there is a
constant $C$ such that for $t \le 1$, on the event $\{ \chi^N_t
(B) \neq 0 \}$,
\begin{align}\label{cutoff2}
\mid \Lambda_I (B)_t \mid \le C (\rho +t), \quad d(I) \le N.
\end{align}

For $N \ge 1$, let us then consider
\[
\mathbf{P}^N_t=\mathbb{E} \left( \chi^N_t (B) \exp \left( \sum_{I,
d(I) \le N } \Lambda_I (B)_t \nabla_{I} \right) \right).
\]

For instance
\[
\mathbf{P}^1_t=\mathbb{E} \left( \chi \left( \frac{\mid B_t
\mid^2}{\rho^2} \right) \exp \left(\sum_{i=1}^d B^i_t \nabla_{i}
\right) \right),
\]
and
\[
\mathbf{P}^2_t=\mathbb{E} \left( \chi^2_t (B) \exp \left(
\sum_{i=0}^d B^i_t \nabla_{i} +\frac{1}{2}\sum_{1\leq i < j \leq
d} \left( \int_0^t B^i_s dB^j_s-B^j_s dB^i_s \right)
[\nabla_{i},\nabla_{j}] \right) \right).
\]

The meaning of this last notation is the following. If $f \in
\Gamma ( \mathbb{M}, \mathcal{E} )$, then $ (\mathbf{P}^N_t
f)(x)=\mathbb{E} (\chi^N_t (B) \Psi (1 ,x))$, where $\Psi (\tau ,
x)$ is the solution of the first order partial differential
equation with random coefficients:
\[
\frac{\partial \Psi}{\partial \tau}(\tau ,x)=\sum_{I, d(I) \le N }
\Lambda_I (B)_t (\nabla_{I} \Psi) (\tau ,x), \quad \Psi
(0,x)=f(x).
\]

Let us consider the following family of norms: If $f \in \Gamma (
\mathbb{M}, \mathcal{E} ) $, for $k \ge 0$,
\[
\parallel f \parallel_k =\sup_{0 \le l \le k} \sup_{0\le i_1, \cdots, i_l \le d}
\sup_{x \in \mathbb{M}} \parallel \nabla_{i_1} \cdots \nabla_{i_l}
f (x) \parallel.
\]

\begin{theorem}\label{Parametrix} Let $N \ge 1$ and $k \ge 0$.
For $f \in \Gamma ( \mathbb{M}, \mathcal{E} ) $,
\[
\parallel \mathbf{P}_t f -\mathbf{P}_t^N f \parallel_k=O\left( t^{\frac{N+1}{2}} \right),
\quad t \rightarrow 0.
\]
\end{theorem}

\begin{proof}

Let us denote
\[
W_t=\sum_{I, d(I) \le N } \Lambda_I (B)_t \nabla_{I}.
\]
First, by expanding out the exponential with Taylor formula we
obtain
\[
\exp \left( W_t \right)f=\left( \sum_{j=0}^N \frac{1}{j!}
W_t^j\right)f +t^{\frac{N+1}{2}} \mathbf{R}^1_N (t),
\]
with
\[
t^{\frac{N+1}{2}} \mathbf{R}^1_N (t)=\frac{1}{N!} \int_0^1
(1-\tau)^N \exp \left( \tau W_t \right) W_t^{N+1} f d\tau.
\]
On the event $\{ \chi^N_t (B) \neq 0 \}$, due to (\ref{cutoff2}),
the coefficients of the first order differential operator $W_t$
are bounded, as well as all their derivatives, by a deterministic
constant. Therefore, on this event,
\[
t^{\frac{N+1}{2}} \parallel  \mathbf{R}^1_N (t)
\parallel_k \le C \left( \sum_{I, d(I) \le N } \mid \Lambda_I (B)_t
\mid \right)^{N+1},
\]
where the constant $C$ only depends on $f$. By using the scaling
property of Brownian motion, we deduce that $\mathbb{E} \left(
\chi^N_t (B) \parallel \mathbf{R}^1_N (t) \parallel_k \right)$ is
bounded when $t \rightarrow 0$. We now observe that, due to
Theorem \ref{Chen-Strichartz brownien}, the rearrangement of terms
in the previous formula gives

\[
\left( \sum_{j=0}^N \frac{1}{j!} W_t^j\right)f=f + \sum_{I, d(I)
\le N}\int_{\Delta^{\mid I \mid} [0,t]}  \circ dB^I \nabla_{i_1}
...\nabla_{i_{\mid I \mid}}f +t^{\frac{N+1}{2}} \mathbf{R}^2_N
(t),
\]
where $\mathbb{E} \left( \parallel \mathbf{R}^2_N (t)
\parallel_k^2 \right)$ is bounded when $t \rightarrow 0$. Indeed,
the remainder is a finite sum of terms
\[
\Lambda_{I_1} (B)_t \cdots \Lambda_{I_j} (B)_t \nabla_{I_1}
\cdots \nabla_{I_j} f, \quad d(I_1)+\cdots +d(I_j) \ge N+1,
\]
and $\Lambda_I (B)_t$ has the same law as $t^{\frac{d(I)}{2}}
\Lambda_I (B)_1$. Therefore
\[
\chi^N_t (B) \exp \left( W_t \right)f=\chi^N_t (B) \left( f +
\sum_{I, d(I) \le N}\int_{\Delta^{\mid I \mid} [0,t]} \circ dB^I
\nabla_{i_1} ...\nabla_{i_{\mid I \mid}}f \right)
+t^{\frac{N+1}{2}} \mathbf{R}^3_N (t),
\]
where $\mathbb{E} \left( \parallel \mathbf{R}^3_N (t)
\parallel_k \right)$ is bounded when $t \rightarrow 0$. Since
iterated integrals have moments of all orders, (\ref{cutoff1}) and
Cauchy-Schwarz inequality show that, in the expectation of the
first term of the right hand side, the cutoff $\chi^N_t (B)$ can
be removed up to a term which is $O(t^m)$ for every $m \ge 0$.
Hence
\[
\mathbf{P}_t^N f=f + \sum_{I, d(I) \le N}\mathbb{E}\left(
\int_{\Delta^{\mid I \mid} [0,t]} \circ dB^I \right) \nabla_{i_1}
...\nabla_{i_{\mid I \mid}}f +t^{\frac{N+1}{2}} \mathbf{R}^4_N
(t),
\]
where $ \parallel \mathbf{R}^4_N (t) \parallel_k $ is bounded when
$t \rightarrow 0$. We now have to compute the expectation of
iterated Stratonovitch integrals. An easy computation shows that
if $\mathcal{I}_n$ is the set of words  with length $n$ obtained
by all the possible concatenations of the words
\[
\{ 0 \}, \{ (i,i) \}, \quad i \in \{1,...,d\},
\]
\begin{enumerate}
\item If $I \notin \mathcal{I}_n$ then
\[
\mathbb{E} \left( \int_{\Delta^n [0,t]}  \circ dB^I \right) =0 ;
\]
\item If $I \in \mathcal{I}_n$ then
\[
\mathbb{E} \left( \int_{\Delta^n [0,t]}  \circ dB^I \right)
=\frac{t^{\frac{n+n(I)}{2}}}{2^{\frac{n-n(I)}{2}}\left(\frac{n+n(I)}{2}
\right) ! },
\]
where $n(I)$ is the number of 0 in $I$ (observe that since $I \in
\mathcal{I}_n$, $n$ and $n(I)$ necessarily have the same parity).
\end{enumerate}
We conclude therefore
\[
\parallel \mathbf{P}^N_t f -\sum_{j  \le \frac{N}{2}} \frac{t^j}{j!} \mathcal{L}^j f\parallel_k =
O\left( t^{\frac{N+1}{2}} \right).
\]
Since it is known that
\[
\parallel \mathbf{P}_t f -\sum_{j  \le \frac{N}{2}} \frac{t^j}{j!} \mathcal{L}^j f \parallel_k =
O\left( t^{\frac{N+1}{2}} \right),
\]
the theorem is proved.

\end{proof}

Let us now fix $x_0 \in \mathbb{M}$ and turn to pointwise
estimates at $x_0$. For these estimates, it is enough to assume
that $\mathcal{L}$ can be written as above in a neighborhood of
$x_0$ only, and this generality will be needed in Section 4. More
precisely, in the remainder of this section, we work under the
following assumptions:
\begin{enumerate}
\item $\mathcal{L}$ is a second order differential operator on
$\mathcal{E}$ such that $-2\mathcal{L}$ is a generalized Laplacian
in the sense of Chapter 2 of \cite{Ber-Ge-Ve}, that is, the
principal symbol of $-2\mathcal{L}$ is given by the Riemannian
metric of $\mathbb{M}$;
\item There exists an open neighborhood $\mathcal{O}$ of $x_0$ on
which
\[
\mathcal{L}=\nabla_{0}+\frac{1}{2} \sum_{i=1}^d \nabla_{i}^2,
\quad \nabla_i=\mathcal{F}_i +\nabla_{V_i}, \quad 0 \le i \le d,
\]
where, as before, the $V_i$'s are smooth vector fields on
$\mathbb{M}$ and the $\mathcal{F}_i$'s are smooth potentials.
\end{enumerate}
The first assumption implies that $\mathcal{L}$ is elliptic and
that $(V_1 ,...,V_d )$ is an orthonormal frame on $\mathcal{O}$.
The operators $\mathbf{P}^N_t$ are defined as before from the
$\nabla_i$'s, and $\mathbf{P}_t=e^{t\mathcal{L}}$. In that case,
$\mathbf{P}_t$ is known to admit a smooth kernel $p_t$ with
respect to the Riemannian measure $dy$ of $\mathbb{M}$: For $f \in
\Gamma ( \mathbb{M}, \mathcal{E} ) $,
\[
(\mathbf{P}_t f)(x_0)=\int_\mathbb{M} p_t (x_0,y)f(y)dy.
\]
Moreover (see Theorem 2.30 in \cite{Ber-Ge-Ve}), for $y$
sufficiently close to $x_0$ and for every $K \ge 0$, $p_t (x_0,y)$
admits a development
\begin{align}\label{devpt}
p_t (x_0,y)=\frac{e^{-\frac{d^2(x_0,y)}{2t}}}{(2\pi t)^{d/2}}
\sum_{k=0}^K \Psi_k(x_0,y) t^k + t^{K-\frac{d}{2}} \mathbf{R}_K
(t,x_0,y) ,
\end{align}
where the remainder term $\mathbf{R}_K (t,x_0,y)$ is bounded when
$t \rightarrow 0$, uniformly in $x_0$ and $y$, $\Psi_k (x_0,
\cdot)$ is a smooth section of $\mathbf{Hom}
(\mathcal{E},\mathcal{E}_{x_0})$ defined around $x_0$ and
$d(\cdot,\cdot)$ is the Riemannian distance.

\begin{theorem}\label{main-estimation}
Let $N \ge 1$. If the constant $\rho$ is small enough, then for
$t>0$ small enough, there exists a smooth section
\[
p^N_t(x_0,\cdot) \in \Gamma( \mathbb{M}, \mathbf{Hom}
(\mathcal{E},\mathcal{E}_{x_0})),
\]
with support included in $\mathcal{O}$, such that for $f \in
\Gamma ( \mathbb{M}, \mathcal{E} ) $,
\[
(\mathbf{P}^N_t f)(x_0)=\int_\mathbb{M} p^N_t (x_0,y)f(y)dy.
\]
Moreover,
\[
p_t (x_0,x_0)=p_t^N (x_0,x_0)+O\left( t^{\frac{N-d}{2}} \right),
\quad t \rightarrow 0.
\]
\end{theorem}

Let us observe that there is a loss of $t^{\frac{1}{2}}$ with
respect to what the weak estimate of Theorem \ref{Parametrix}
could suggest. The truncation of the Lie series neglects terms of
order $t^{\frac{N+1}{2}}$ in the position of a flow that lives at
the scale $t^{\frac{1}{2}}$, which is a relative error of order
$t^{\frac{N}{2}}$, and the kernels are of order
$t^{-\frac{d}{2}}$. For another approach to the diagonal
asymptotics of heat kernels by stochastic Taylor expansions, that
is based on the Malliavin calculus (see \cite{Mal2}, \cite{Nu}),
we refer to \cite{Ben}.

\begin{proof}

The proof is not simple. We shall proceed in several steps. In a
first step, we use the Brownian bridges $B^{t,b}$ to show the
existence of the kernel $p^N_t (x_0,\cdot)$ and to give an
explicit expression of it. In a second step, we show that, at the
scale $\sqrt{t}$, the kernels $p^N_t (x_0,\cdot)$ and
$p_t(x_0,\cdot)$ both admit a development in powers of $\sqrt{t}$.
And finally, we shall prove the required estimate by comparing
these two developments, thanks to a refinement of Theorem
\ref{Parametrix}.

\

In the whole proof, we use the normal coordinates
$y=(y^1,...,y^d)$ centered at $x_0$ that are associated with the
orthonormal basis $(V_1(x_0),...,V_d(x_0))$. Shrinking
$\mathcal{O}$ if needed, we therefore identify $\mathcal{O}$ with
a ball of $\mathbb{R}^d$ centered at $x_0=0$, in such a way that
\[
V_i (x_0)=\frac{\partial}{\partial y^i}, \quad d(x_0,y)=\mid y
\mid, \quad dy=\theta (y) dy^1 \cdots dy^d,
\]
where $\theta$ is a smooth function such that $\theta (0)=1$. We
also fix a trivialization of $\mathcal{E}$ over $\mathcal{O}$:
The sections of $\mathcal{E}$ over $\mathcal{O}$ are then functions
with values in $\mathcal{E}_{x_0}$, whereas $p_t (x_0,y)$ and
$\Psi_k (x_0,y)$ are elements of $\mathbf{End}
(\mathcal{E}_{x_0})$. In what follows, $C$ and $L$ denote
deterministic constants whose value may change from line to line.

\

\textbf{First step:}

\

If $\lambda=(\lambda_I)_{d(I) \le N}$ is a family of real numbers
that are small enough, we denote
\[
\phi (\lambda)=\exp \left( \sum_{I, d(I) \le N } \lambda_I V_{I}
\right) (x_0) \in \mathcal{O}.
\]
If $f$ is a smooth function on $\mathbb{M}$, we denote by
$\mathcal{M}_f$ the operator on $\Gamma (\mathbb{M} ,
\mathcal{E})$ that acts by multiplication by $f$. Due to the
Leibniz rule for connections, we have for any word $I$:
\[
[\nabla_I, \mathcal{M}_f] =\mathcal{M}_{V_I f}.
\]
Consequently,
\[
\left[ \sum_{I, d(I) \le N } \lambda_I \nabla_{I} ,\mathcal{M}_f
\right]= \mathcal{M}_{\sum_{I, d(I) \le N } \lambda_I V_{I} f}.
\]
The above commutation property implies the following one:
\[
\exp \left( \sum_{I, d(I) \le N } \lambda_I \nabla_{I} \right)
\mathcal{M}_f=\mathcal{M}_{\exp \left( \sum_{I, d(I) \le N }
\lambda_I V_{I} \right) f}\exp \left( \sum_{I, d(I) \le N }
\lambda_I \nabla_{I} \right).
\]
Therefore, there exists $\mathcal{H} (\lambda) \in \mathbf{End}
(\mathcal{E}_{x_0})$ such that for every $\eta \in \Gamma
(\mathbb{M} , \mathcal{E})$,
\[
\left[\exp \left( \sum_{I, d(I) \le N } \lambda_I \nabla_{I}
\right)\eta \right] (x_0)=\mathcal{H} (\lambda) \eta (\phi
(\lambda)).
\]
Since the solutions of ordinary differential equations smoothly
depend on parameters, $\phi$ and $\mathcal{H}$ are smooth
functions of $\lambda$. We have $\mathcal{H}(0)=\mathbf{Id}$ and
we can write
\begin{align}\label{dev-phi}
\phi (\lambda)=\sum_{I, d(I) \le N} \lambda_I V_I (x_0) +\sum_{I,
I'} \lambda_I \lambda_{I'} \phi_{I,I'} (\lambda),
\end{align}
where the $\phi_{I,I'}$'s are smooth.

\

Let us now come back to the Brownian bridges $B^{t,b}$ and denote
\[
\lambda_t (b)=\left( \Lambda_I (B^{t,b})_t \right)_{d(I) \le N},
\quad \Phi_t (b)=\phi \left( \lambda_t (b) \right), \quad
\mathcal{H}_t (b)= \mathcal{H} \left( \lambda_t (b) \right), \quad
\chi_t (b)=\chi^N_t (B^{t,b}),
\]
that is
\[
\chi_t (b)=\chi \left( \frac{1}{\rho^2} \left( \mid b \mid^2
+\sum_{J} \frac{\Lambda_J (B^{t,0})_t^2}{t^{2n(J)}}
\right)\right).
\]
These are random functions of $b$ that only depend on the loop
$B^{t,0}$, and due to (\ref{cutoff2}), $\Phi_t (b)$ and
$\mathcal{H}_t (b)$ are well defined on the support of $\chi_t$
if $\rho$ and $t$ are small enough. Since $B^{t,B_t}=B$, the
definition of $\mathbf{P}^N_t$ reads
\[
(\mathbf{P}^N_t f)(x_0)=\mathbb{E} \left( \chi_t (B_t)
\mathcal{H}_t (B_t) f (\Phi_t (B_t)) \right),
\]
and since $B_t$ is a Gaussian random variable which is
independent of $B^{t,0}$,
\begin{align}\label{disintegration}
(\mathbf{P}^N_t f)(x_0)=\mathbb{E} \left( \int_{\mathbb{R}^d}
\chi_t (b) \mathcal{H}_t (b) f (\Phi_t (b)) \frac{e^{-\frac{\mid b
\mid^2}{2t}}}{(2\pi t)^{d/2}} db\right).
\end{align}
We now claim that if $\rho$ is small enough, then for $t \le
\rho$, on the event
\[
\mathcal{A}_t= \left\{ \sum_{J} \frac{\Lambda_J
(B^{t,0})_t^2}{t^{2n(J)}} \le 4 \rho^2 \right\}
\]
outside of which $\chi_t$ identically vanishes, the map $b
\rightarrow \Phi_t (b)$ is a diffeomorphism from the ball $\{ \mid
b \mid < 3\rho \}$ onto an open subset of $\mathcal{O}$ and
satisfies
\begin{align}\label{jacobian}
\parallel D_b \Phi_t (b) -\mathbf{Id} \parallel \le \frac{1}{2}.
\end{align}
Indeed, as for (\ref{cutoff2}), we see from (\ref{polynomial})
that on $\mathcal{A}_t$, for $\mid b \mid < 3\rho$ and $\mid I
\mid \ge 2$,
\[
\mid \Lambda_I (B^{t,b})_t \mid \le C (\rho +t), \quad \left|
\frac{\partial }{\partial b_j} \Lambda_I (B^{t,b})_t \right| \le C
(\rho +t),
\]
whereas $\Lambda_{(i)} (B^{t,b})_t=b_i$ and $\Lambda_{(0)}
(B^{t,b})_t=t$. Since $V_i (x_0)=\frac{\partial}{\partial y^i}$,
we deduce from (\ref{dev-phi}) that
\[
\Phi_t (b)=b +\sum_{I, 2 \le d(I) \le N} \Lambda_I
(B^{t,b})_t V_I (x_0) +\sum_{I, I'} \Lambda_I (B^{t,b})_t
\Lambda_{I'} (B^{t,b})_t \phi_{I,I'} \left( \lambda_t (b) \right)
\]
satisfies $\parallel D_b \Phi_t (b) -\mathbf{Id} \parallel \le C
(\rho+t)$, which gives (\ref{jacobian}). As a consequence, $b
\rightarrow \Phi_t (b)-b$ is a contraction on the ball, so that
$\Phi_t$ is injective, and it is therefore a diffeomorphism onto
its image.

\

Since $\chi_t (b)=0$ for $\mid b \mid \ge 2 \rho$, we can
therefore perform in (\ref{disintegration}) the change of variable
$y=\Phi_t (b)$. We obtain
\[
(\mathbf{P}^N_t f)(x_0)=\int_\mathbb{M} p^N_t (x_0,y)f(y)dy,
\]
with
\begin{align}\label{kernel-formula}
p^N_t (x_0,y)=\frac{1}{\theta (y) (2\pi t)^{d/2}} \mathbb{E}
\left( \left[ \chi_t (b) \frac{e^{-\frac{\mid b \mid^2}{2t}}}{\mid
\det D_b \Phi_t (b) \mid} \mathcal{H}_t (b)\right]_{b=\Phi_t^{-1}
(y)} \right),
\end{align}
where it is understood that the quantity inside the expectation is
$0$ outside of $\mathcal{A}_t$ or if $y$ is not in the image of
$\Phi_t$. On $\mathcal{A}_t$ and on the ball $\{ \mid b \mid <
3\rho \}$, the derivatives at any order of $\Phi_t$ and
$\mathcal{H}_t$ are bounded by deterministic constants and $\mid
\Phi_t (b) \mid \le C (\rho +t)$. Therefore $p^N_t (x_0,\cdot)$ is
smooth and, if $\rho$ and $t$ are small enough, its support is
included in $\mathcal{O}$.

\

\textbf{Second step:}

\

Let us now study $p^N_t (x_0,\cdot)$ when $t \rightarrow 0$. We
denote $\varepsilon=\sqrt{t}$. By the scaling property of Brownian
motion, the family $\left( \Lambda_I (B^{t,\varepsilon b})_t
\right)_{I,b}$ has the same law as the family $\left(
\varepsilon^{d(I)} \Lambda_I (B^{1,b})_1 \right)_{I,b}$.
Therefore, if we denote
\[
\bar{\lambda}_\varepsilon (b)=\left( \varepsilon^{d(I)} \Lambda_I
(B^{1,b})_1 \right)_{d(I) \le N}, \quad F_\varepsilon
(b)=\frac{1}{\varepsilon} \phi \left( \bar{\lambda}_\varepsilon
(b) \right), \quad \bar{\mathcal{H}}_\varepsilon (b)= \mathcal{H}
\left( \bar{\lambda}_\varepsilon (b) \right),
\]
and
\[
\bar{\chi}_\varepsilon (b)=\chi \left( \frac{1}{\rho^2} \left(
\varepsilon^2 \mid b \mid^2 +\sum_{J} \varepsilon^{2(\mid J \mid
-n(J))} \Lambda_J (B^{1,0})_1^2 \right)\right),
\]
then the random functions $b \rightarrow \left(
\frac{1}{\varepsilon} \Phi_t (\varepsilon b) , \mathcal{H}_t
(\varepsilon b), \chi_t (\varepsilon b) \right)$ and $b
\rightarrow \left( F_\varepsilon (b) ,
\bar{\mathcal{H}}_\varepsilon (b), \bar{\chi}_\varepsilon (b)
\right)$ have the same law. Consequently, from
(\ref{disintegration}),
\begin{align}\label{scaled-weak}
(\mathbf{P}^N_t f)(x_0)=\mathbb{E} \left( \bar{\chi}_\varepsilon
(B_1) \bar{\mathcal{H}}_\varepsilon (B_1) f (\varepsilon
F_\varepsilon (B_1)) \right),
\end{align}
and from (\ref{kernel-formula}),
\begin{align}\label{scaled-kernel}
\varepsilon^d \theta (\varepsilon z) p^N_t (x_0,\varepsilon z)=
\mathbb{E} \left( \left[ \bar{\chi}_\varepsilon (b) \frac{\gamma
(b)}{\mid \det D_b F_\varepsilon (b) \mid}
\bar{\mathcal{H}}_\varepsilon (b)\right]_{b=F_\varepsilon^{-1}
(z)} \right), \quad \gamma (b)=\frac{e^{-\frac{\mid b
\mid^2}{2}}}{(2\pi )^{d/2}}.
\end{align}
The interest of this scaling is that everything is now smooth with
respect to $\varepsilon$, up to $\varepsilon=0$. Indeed, from
(\ref{dev-phi}),
\[
F_\varepsilon (b)=b +\sum_{I, 2 \le d(I) \le N} \varepsilon^{d(I)-1}
\Lambda_I (B^{1,b})_1 V_I (x_0) +\sum_{I, I'}
\varepsilon^{d(I)+d(I')-1} \Lambda_I (B^{1,b})_1 \Lambda_{I'}
(B^{1,b})_1 \phi_{I,I'} \left( \bar{\lambda}_\varepsilon (b)
\right),
\]
and from (\ref{polynomial}), the $\Lambda_I (B^{1,b})_1$'s are
polynomials in $b$ whose coefficients are the $\Lambda_J
(B^{1,0})_1$'s. We deduce that, where they are defined,
$F_\varepsilon (b)$, $\bar{\mathcal{H}}_\varepsilon (b)$ and
$\bar{\chi}_\varepsilon (b)$ are smooth functions of
$(\varepsilon,b)$, up to $\varepsilon=0$, with
\[
F_0 (b)=b, \quad \bar{\mathcal{H}}_0 (b)=\mathbf{Id}, \quad
\bar{\chi}_0 (b)=1.
\]
Moreover:
\begin{enumerate}
\item On the support of $\bar{\chi}_\varepsilon$, and more
generally for $\varepsilon \mid b \mid <3 \rho$ on the event
\[
\bar{\mathcal{A}}_\varepsilon= \left\{ \sum_{J}
\varepsilon^{2(\mid J \mid -n(J))} \Lambda_J (B^{1,0})_1^2 \le 4
\rho^2 \right\}
\]
that corresponds to $\mathcal{A}_t$, the derivatives at any order
with respect to $(\varepsilon,b)$ of these three functions are
bounded by $C(1+\mid b \mid +\Xi)^L$, where
\[
\Xi=\sum_{J, \mid J \mid \le N} \mid \Lambda_J (B^{1,0})_1 \mid
\]
has moments of all orders;
\item The derivatives at any order greater than or equal to $1$
with respect to $\varepsilon$ of $\bar{\chi}_\varepsilon (b)$
vanish at $\varepsilon=0$, and
those of $F_\varepsilon (b)$ and $\bar{\mathcal{H}}_\varepsilon
(b)$ at $\varepsilon=0$ are polynomials in $b$ and in the
$\Lambda_J (B^{1,0})_1$'s;
\item From (\ref{jacobian}), on $\bar{\mathcal{A}}_\varepsilon$
and for $\varepsilon \mid b \mid <3 \rho$, $\parallel D_b
F_\varepsilon (b) -\mathbf{Id} \parallel \le \frac{1}{2}$, so that
\begin{align}\label{inverse-bound}
\mid F_\varepsilon^{-1} (z) \mid \le 2 \mid z -F_\varepsilon (0)
\mid, \quad \mid F_\varepsilon (0) \mid \le C \varepsilon (1+\Xi).
\end{align}
\end{enumerate}
We now draw two consequences of these properties.

\

The first one is a weak development. Let $g: \mathbb{R}^d
\rightarrow \mathbb{R}$ be a smooth function whose derivatives
have a polynomial growth. From Taylor formula with respect to
$\varepsilon$ and from the first two properties, we get
\begin{align}\label{weak-N}
\mathbb{E} \left( \bar{\chi}_\varepsilon (B_1)
\bar{\mathcal{H}}_\varepsilon (B_1) g ( F_\varepsilon (B_1))
\right)=\sum_{n=0}^{N-1} \varepsilon^n \int_{\mathbb{R}^d} \gamma
(z) \tilde{\pi}_n (z) g(z) dz +O\left(\varepsilon^N \right),
\end{align}
where the $\tilde{\pi}_n$'s are polynomials with values in
$\mathbf{End} (\mathcal{E}_{x_0})$ that do not depend on $g$.
Indeed, the derivative of order $n$ at $\varepsilon=0$ of
$\bar{\mathcal{H}}_\varepsilon (b) g ( F_\varepsilon (b))$ is a
finite sum of terms $c(b) \partial^\alpha g (b)$, where $c(b)$ is
a polynomial in $b$ and in the $\Lambda_J (B^{1,0})_1$'s. Since
$B^{1,0}$ is independent of $B_1$, we can first take the
expectation with respect to $B^{1,0}$, and then integrate by parts
with respect to the Gaussian measure $\gamma (b) db$. This leads
to (\ref{weak-N}).

\

The second consequence is a pointwise development. By the implicit
function theorem and (\ref{inverse-bound}), $F_\varepsilon^{-1}
(z)$ is a smooth function of $(\varepsilon,z)$, up to
$\varepsilon=0$, whose derivatives are bounded by $C(1+\mid z \mid
+\Xi)^L$, and $F_0^{-1} (z)=z$. The quantity inside the
expectation in (\ref{scaled-kernel}) is therefore a smooth
function of $(\varepsilon,z)$ whose derivatives satisfy the same
type of bound (let us recall that $\bar{\chi}_\varepsilon (b)=0$
for $\varepsilon \mid b \mid \ge 2 \rho$, whereas $F_\varepsilon$
is a diffeomorphism on $\{ \varepsilon \mid b \mid <3\rho \}$).
From Taylor formula we deduce that, uniformly for $z$ in compact
sets,
\begin{align}\label{pointwise-N}
\varepsilon^d \theta (\varepsilon z) p^N_t (x_0,\varepsilon z)=
\sum_{n=0}^{N-1} \varepsilon^n \tilde{a}_n (z)
+O\left(\varepsilon^N \right),
\end{align}
where the $\tilde{a}_n$'s are continuous. But if $g$ is smooth and
compactly supported, by applying (\ref{scaled-weak}) to $f(y)=g
\left( \frac{y}{\varepsilon} \right) u$ for every $u \in
\mathcal{E}_{x_0}$, we see that the left hand side of
(\ref{weak-N}) is also equal to the endomorphism
\[
\int_{\mathbb{R}^d} \varepsilon^d \theta (\varepsilon z) p^N_t
(x_0,\varepsilon z) g(z) dz.
\]
By comparing (\ref{weak-N}) and (\ref{pointwise-N}) we obtain
therefore $\tilde{a}_n=\gamma \tilde{\pi}_n$. In particular, at
$z=0$,
\begin{align}\label{diag-N}
t^{\frac{d}{2}} p^N_t (x_0,x_0)=\frac{1}{(2\pi )^{d/2}}
\sum_{n=0}^{N-1} t^{\frac{n}{2}} \tilde{\pi}_n (0)
+O\left(t^{\frac{N}{2}} \right).
\end{align}

\

Let us finally turn to the kernel $p_t$. Since $d(x_0,y)=\mid y
\mid$, the development (\ref{devpt}) reads
\[
\varepsilon^d \theta (\varepsilon z) p_t (x_0,\varepsilon
z)=\gamma (z) \theta (\varepsilon z) \sum_{k=0}^K
\varepsilon^{2k} \Psi_k(x_0,\varepsilon z) + \varepsilon^{2K}
\theta (\varepsilon z) \mathbf{R}_K (t,x_0,\varepsilon z).
\]
By expanding the smooth functions $\theta \Psi_k (x_0,\cdot)$ at
$0$ with Taylor formula, we deduce that for $2K \ge N$ and
$\varepsilon z \in \mathcal{O}$,
\begin{align}\label{true-N}
\varepsilon^d \theta (\varepsilon z) p_t (x_0,\varepsilon
z)=\gamma (z) \left( \sum_{n=0}^{N-1} \varepsilon^n \pi_n (z)
+\varepsilon^N \mathbf{r}_N (\varepsilon,z) \right)+
\varepsilon^{2K} \mathbf{r}'_K (\varepsilon,z),
\end{align}
where $\mid \mathbf{r}_N (\varepsilon,z) \mid \le C (1+\mid z
\mid)^N$, $\mathbf{r}'_K$ is bounded, and $\pi_n$ is the
polynomial
\[
\pi_n (z)=\sum_{2k+r=n} \left[ \theta \Psi_k (x_0,\cdot)
\right]_r (z),
\]
$[ \cdot ]_r$ denoting the homogeneous part of degree $r$ in the
Taylor development at $0$. In particular, since $\theta (0)=1$,
\begin{align}\label{diag-true}
t^{\frac{d}{2}} p_t (x_0,x_0)=\frac{1}{(2\pi )^{d/2}}
\sum_{n=0}^{N-1} t^{\frac{n}{2}} \pi_n (0) +O\left(t^{\frac{N}{2}}
\right).
\end{align}

\

\textbf{Final step:}

\

From (\ref{diag-N}) and (\ref{diag-true}), the theorem is proved
if we show that
\[
\pi_n=\tilde{\pi}_n, \quad 0 \le n \le N-1.
\]
To this end, let $g$ be a homogeneous polynomial function on
$\mathbb{R}^d$ with degree $p \ge 0$, let $u \in
\mathcal{E}_{x_0}$, and let $\varphi$ be a smooth function with
compact support included in $\mathcal{O}$ that is equal to $1$ in
a neighborhood of $x_0$. In our coordinates, $f=\varphi g u$ is a
smooth section of $\mathcal{E}$ over $\mathbb{M}$ that vanishes at
order $p$ at $x_0$, that is, $f$ and its derivatives of order less
than $p$ vanish at $x_0$. We claim that for such a section,
\begin{align}\label{refined}
(\mathbf{P}_t f)(x_0) -(\mathbf{P}_t^N f) (x_0)=O\left(
t^{\frac{N+p}{2}} \right), \quad t \rightarrow 0.
\end{align}
To prove this claim, we come back to the proof of Theorem
\ref{Parametrix}. With the same notations, Taylor formula at the
order $N+p$ gives
\[
\exp \left( W_t \right)f=\left( \sum_{j=0}^{N+p} \frac{1}{j!}
W_t^j\right)f +t^{\frac{N+p+1}{2}} \mathbf{R} (t),
\]
where $\mathbb{E} \left( \chi^N_t (B) \parallel \mathbf{R} (t)
\parallel_0 \right)$ is bounded. Let us give to $X_0$ the weight
$2$ and to $X_1,...,X_d$ the weight $1$, so that $X_I$ has the
weight $d(I)$. In $\mathbb{R} [[X_0,...,X_d]]$, the components
with weight less than $N+p$ of
\[
\sum_{j=0}^{N+p} \frac{1}{j!} \left( \sum_{I, d(I) \le N}
\Lambda_I (B)_t X_I \right)^j \quad \text{and} \quad \exp \left(
\sum_{k \geq 1} \sum_{I \in \{0,1,...,d\}^k}\Lambda_I (B)_t X_I
\right)
\]
only differ by products
\[
\Lambda_{I_1} (B)_t \cdots \Lambda_{I_j} (B)_t X_{I_1} \cdots
X_{I_j}
\]
such that $d(I_l) \ge N+1$ for at least one of the words. The
weight of such a product is at least $N+j$, and it is less than
$N+p$ only if $j <p$. But since the $\nabla_I$'s are first order
differential operators and $f$ vanishes at order $p$ at $x_0$,
\[
(\nabla_{I_1} \cdots \nabla_{I_j} f)(x_0)=0, \quad j <p.
\]
Therefore, from Theorem \ref{Chen-Strichartz brownien},
\[
\left[ \left( \sum_{j=0}^{N+p} \frac{1}{j!} W_t^j\right)f \right]
(x_0)=f(x_0) + \sum_{I, d(I) < N+p}\int_{\Delta^{\mid I \mid}
[0,t]}  \circ dB^I (\nabla_{i_1} ...\nabla_{i_{\mid I \mid}}f)
(x_0) +t^{\frac{N+p}{2}} \mathbf{R}' (t),
\]
where $\mathbb{E} \left( \mid \mathbf{R}' (t) \mid^2 \right)$ is
bounded. By taking the expectation as in the proof of Theorem
\ref{Parametrix}, and by using the fact that
$\mathcal{L}=\nabla_{0}+\frac{1}{2} \sum_{i=1}^d \nabla_{i}^2$ on
$\mathcal{O}$, we obtain
\[
(\mathbf{P}_t^N f) (x_0)= \sum_{j < \frac{N+p}{2}} \frac{t^j}{j!}
(\mathcal{L}^j f) (x_0) +O\left( t^{\frac{N+p}{2}} \right).
\]
Since the same development holds for $(\mathbf{P}_t f)(x_0)$, the
claim (\ref{refined}) is proved.

\

We can now conclude. On one hand, for $\rho$ and $t$ small enough,
$\varphi=1$ on the support of $p^N_t (x_0,\cdot)$. Since $g$ is
homogeneous, (\ref{scaled-weak}) and (\ref{weak-N}) give
therefore
\[
(\mathbf{P}_t^N f) (x_0)=\varepsilon^p \text{ } \mathbb{E} \left(
\bar{\chi}_\varepsilon (B_1) \bar{\mathcal{H}}_\varepsilon (B_1) g
( F_\varepsilon (B_1)) \right) u =\varepsilon^p \left(
\sum_{n=0}^{N-1} \varepsilon^n \int_{\mathbb{R}^d} \gamma (z)
\tilde{\pi}_n (z) g(z) u dz +O\left(\varepsilon^N \right) \right).
\]
On the other hand,
\[
(\mathbf{P}_t f) (x_0)=\int_{\mathcal{O}} p_t (x_0,y) f(y)
dy=\varepsilon^p \int_{\mathbb{R}^d} \varepsilon^d \theta
(\varepsilon z) p_t (x_0,\varepsilon z) \varphi (\varepsilon z)
g(z) u dz,
\]
and from (\ref{true-N}), used with $2K \ge N+d+p$,
\[
(\mathbf{P}_t f) (x_0)=\varepsilon^p \left( \sum_{n=0}^{N-1}
\varepsilon^n \int_{\mathbb{R}^d} \gamma (z) \pi_n (z) g(z) u dz
+O\left(\varepsilon^N \right) \right)
\]
(in the Gaussian integrals, the function $\varphi (\varepsilon
z)$ can be replaced by $1$ up to a term which is exponentially
small when $\varepsilon \rightarrow 0$).
From (\ref{refined}), we deduce that for $0 \le n \le N-1$,
\[
\int_{\mathbb{R}^d} \gamma (z) \left( \pi_n (z) -\tilde{\pi}_n (z)
\right) g(z) u dz=0.
\]
This holds for every homogeneous polynomial function $g$, and thus
for every polynomial function $g$. Since $\pi_n -\tilde{\pi}_n$ is
itself a polynomial, this implies $\pi_n=\tilde{\pi}_n$.
\end{proof}

If $I\in \{0,1,...,d\}^k$ is a word, we denote
\[
\mathcal{F}_I=\nabla_I-\nabla_{V_I} \in \Gamma
(\mathbb{M},\mathbf{End}(\mathcal{E})).
\]
The value of $p^N_t (x_0,x_0)$ can be made explicit in terms of
these potentials.

\begin{corollary}\label{corollary-diag}
Let $N \ge 1$. If the constant $\rho$ is small enough, when $t
\rightarrow 0$,
\[
p_t (x_0,x_0)=\frac{1}{(2\pi t)^{d/2}} \mathbb{E} \left( \chi^N_t
(B^{t,b^*_t}) \frac{e^{-\frac{\mid b^*_t \mid^2}{2t}}}{\mid \det
\mathbf{J}_t \mid} \exp \left( \sum_{I, d(I) \le N } \Lambda_I
(B^{t,b^*_t})_t \mathcal{F}_{I}(x_0) \right) \right) + O\left(
t^{\frac{N-d}{2}}\right),
\]
where $b^*_t$ is the solution $b$, $\mid b \mid <3\rho$, of the
equation
\[
\sum_{I, d(I) \le N } \Lambda_I (B^{t,b})_t V_{I}(x_0)=0,
\]
which is unique as soon as $b \rightarrow \chi^N_t (B^{t,b})$ does
not identically vanish, where $\mathbf{J}_t$ is the differential
at $b^*_t$ of the map
\[
b \rightarrow \exp \left( \sum_{I, d(I) \le N } \Lambda_I
(B^{t,b})_t V_{I} \right)(x_0)
\]
computed in normal coordinates at $x_0$, and where the quantity
inside the expectation is understood to be $0$ if $b \rightarrow
\chi^N_t (B^{t,b})$ identically vanishes or if there is no such
$b^*_t$.
\end{corollary}

\begin{proof}
With the notations of the proof of Theorem \ref{main-estimation},
this theorem and the formula (\ref{kernel-formula}) taken at
$y=x_0$ give
\[
p_t (x_0,x_0)=\frac{1}{(2\pi t)^{d/2}} \mathbb{E} \left( \left[
\chi_t (b) \frac{e^{-\frac{\mid b \mid^2}{2t}}}{\mid \det D_b
\Phi_t (b) \mid} \mathcal{H}_t (b)\right]_{b=\Phi_t^{-1} (x_0)}
\right)+ O\left( t^{\frac{N-d}{2}}\right).
\]
We therefore have to understand the equation $\Phi_t (b)=x_0$ and
to compute $\mathcal{H}_t (b)$ when $b$ is a solution of it. So,
let $\lambda=(\lambda_I)_{d(I) \le N}$ be such that $\mid
\lambda_I \mid \le C(\rho +t)$ and let $V=\sum_{I, d(I) \le N}
\lambda_I V_I$. We claim that if $\rho$ and $t$ are small enough,
then
\[
\exp (V) (x_0)=x_0 \quad \text{if and only if} \quad V(x_0)=0.
\]
Indeed, in the normal coordinates, let $x_s=\exp (sV)(x_0)$ and
let $\kappa$ be the supremum on $\mathcal{O}$ of the norm of the
differential of $V$, so that $\kappa \le C (\rho+t)$. Since $\mid
V(x_s) \mid \le \mid V(x_0) \mid +\kappa \mid x_s -x_0 \mid$,
Gronwall's lemma gives $\mid x_s -x_0 \mid \le e^{\kappa} \mid
V(x_0) \mid$ for $0 \le s \le 1$. Therefore, if $x_1=x_0$,
\[
\mid V(x_0) \mid =\left| \int_0^1 (V(x_0)-V(x_s)) ds \right| \le
\kappa e^{\kappa} \mid V(x_0) \mid,
\]
which implies $V(x_0)=0$ as soon as $\kappa e^{\kappa}<1$. The
converse is obvious. Consequently, $\Phi_t^{-1} (x_0)=b^*_t$ and
$D_b \Phi_t (b^*_t)=\mathbf{J}_t$. Finally, if $V(x_0)=0$, then
$x_0$ is a fixed point of the flow of $V$, so that for $\eta \in
\Gamma (\mathbb{M} , \mathcal{E})$, the function
\[
\tau \rightarrow \left[\exp \left( \tau \sum_{I, d(I) \le N }
\lambda_I \nabla_{I} \right)\eta \right] (x_0)
\]
solves a linear differential equation in $\mathcal{E}_{x_0}$, and
\[
\left[\exp \left( \sum_{I, d(I) \le N } \lambda_I \nabla_{I}
\right)\eta \right] (x_0)=\exp \left( \sum_{I, d(I) \le N }
\lambda_I \mathcal{F}_{I}(x_0) \right) \eta (x_0).
\]
That is, $\mathcal{H} (\lambda)=\exp \left( \sum_{I, d(I) \le N }
\lambda_I \mathcal{F}_{I}(x_0) \right)$, and the corollary is
proved.
\end{proof}

\section{The local index theorem for the Dirac operator on the
spin bundle}

The Atiyah-Singer index theorem for the Dirac operator on the spin
bundle as proved in \cite{At-Si}, is the following:

\begin{theorem}
Let $\mathbb{M}$ be a compact, $d$-dimensional spin manifold, with
$d$ even. Let $\mathbf{D}$ be the Dirac operator on the spin
bundle $\mathcal{S}=\mathcal{S}^+ \oplus \mathcal{S}^-$ of
$\mathbb{M}$ and let $\mathbf{ind} (\mathbf{D} )$ be the index of
its restriction $\mathbf{D}^+: \Gamma (\mathbb{M}, \mathcal{S}^+)
\rightarrow \Gamma (\mathbb{M}, \mathcal{S}^-)$. Then
\[
\mathbf{ind} (\mathbf{D} )=\left( \frac{1}{2i\pi}
\right)^{\frac{d}{2}} \int_\mathbb{M} [A(\mathbb{M})]_d,
\]
where $[A(\mathbb{M})]_d$ is the form of degree $d$ on $\mathbb{M}$
obtained by taking the $d$-form piece of the $A$-genus
\[
A(\mathbb{M})=\det \left( \frac{ \Omega}{2 \sinh \frac{1}{2}
\Omega} \right)^{\frac{1}{2}},
\]
and $\Omega$ is the Riemannian curvature form defined in local
orthonormal frame $e_i$ with dual frame $e^*_i$ by
\[
\Omega=\frac{1}{2} \sum_{1 \le i,j \le d} R(e_i,e_j) e_i^* \wedge
e^*_j,
\]
with $R$, Riemannian curvature.
\end{theorem}

Before we turn to the proof. Let us first recall some linear
algebra constructions as can be found in Chapter 3 of
\cite{Ber-Ge-Ve}.

\

Let $V$ be an oriented $d$ dimensional Euclidean space. We assume
that the dimension $d$ is even. The Clifford algebra
$\mathbf{Cl}(V)$ over $V$ is the algebra
\[
\mathbf{T} (V)=\mathbb{R} \oplus V \oplus (V \otimes V) \oplus
\cdots
\]
quotiented by the relations
\begin{align}\label{relation-clifford}
u \otimes v +v \otimes u +2 \langle u,v \rangle 1 =0.
\end{align}
Let $e_1,...,e_d$ be an oriented orthonormal basis of $V$. The family
\[
e_{i_1} ... e_{i_k}, \quad 0 \le k \le d, \quad 1 \le i_1 <...<i_k
\le d,
\]
forms a basis of $\mathbf{Cl}(V)$ that is therefore of dimension
$2^d$. In $\mathbf{T} (V)$ we can distinguish elements that are
even from elements that are odd. This leads to a decomposition:
\[
\mathbf{Cl}(V)=\mathbf{Cl}^-(V) \oplus \mathbf{Cl}^+(V),
\]
with $V \subset \mathbf{Cl}^-(V)$.

A Clifford module is a vector space $E$ over $\mathbb{R}$ (or
$\mathbb{C}$) that is also a $\mathbf{Cl}(V)$-module and that
admits a direct sum decomposition
\[
E=E^- \oplus E^+
\]
with
\[
\mathbf{Cl}^+(V) \cdot E^\pm \subset E^\pm, \quad \mathbf{Cl}^-(V)
\cdot E^\pm \subset E^\mp.
\]
Up to isomorphism, there is a unique irreducible complex module $S$
for $\mathbb{C} \otimes \mathbf{Cl}(V)$, called the spinor module
over $V$. Its Clifford action $c$ induces the isomorphism
\[
\mathbb{C} \otimes \mathbf{Cl}(V) \simeq \mathbf{End}_{\mathbb{C}} (S).
\]
In particular $\dim_{\mathbb{C}} S=2^{\frac{d}{2}}$. We take $S^\pm$
to be the $\pm 1$ eigenspaces of the chirality operator
\[
\Gamma=i^{\frac{d}{2}} c(e_1)\cdots c(e_d).
\]
The supertrace is then given by
\[
\mathbf{Str} \text{ } a=\mathbf{Tr} \left( \Gamma c(a) \right),
\quad a \in \mathbf{Cl}(V).
\]
If
\[
a =\sum_{k=0}^d \sum_{1 \le i_1 <...<i_k \le d } a_{i_1 \cdots
i_k} e_{i_1} ... e_{i_k},
\]
then we have
\begin{align}\label{supertraceformula-clifford}
\mathbf{Str} \text{ } a=\left( \frac{2}{i} \right)^{\frac{d}{2}}
a_{1 \cdots d}.
\end{align}

If $\psi \in \mathfrak{so} (V)$, that is if $\psi : V \rightarrow
V$ is a skew-symmetric map, we define
\[
D\psi=\frac{1}{2}\sum_{1 \le i<j \le d} \langle \psi (e_i), e_j
\rangle e_i e_j \in \mathbf{Cl}(V),
\]
and observe that $D[\psi_1,\psi_2]=[D\psi_1,D\psi_2]$. The set
$\mathbf{Cl}^2(V)=D \mathfrak{so} (V)$ is therefore a Lie algebra.
The Lie group $\mathbf{Spin} (V)$ is the group obtained by
exponentiating $\mathbf{Cl}^2(V)$ inside the Clifford algebra
$\mathbf{Cl}(V)$; It is a two-fold covering of the special
orthogonal group $\mathbf{SO}(V)$. It can also be described as the
set of $a \in \mathbf{Cl}(V)$ such that:
\[
a=v_1 ... v_{2k}, \quad 1 \le k \le \frac{d}{2},\quad v_i \in V,
\quad \parallel v_i \parallel=1.
\]

We now come back to differential geometry and carry the above
construction on the cotangent spaces of a spin manifold.

\

So, let $\mathbb{M}$ be a compact $d$-dimensional, Riemannian and
oriented manifold. We assume that $d$ is even. We furthermore
assume that $\mathbb{M}$ admits a spin structure: That is, there
exists a principal bundle on $\mathbb{M}$ with structure group
$\mathbf{Spin}(\mathbb{R}^d)$ such that the bundle charts are
compatible with the two-fold covering
$\mathbf{Spin}(\mathbb{R}^d) \rightarrow
\mathbf{SO}(\mathbb{R}^d)$. This bundle will be denoted
$\mathcal{SP} (\mathbb{M})$ and $\pi$ will denote the canonical
surjection. The spin bundle $\mathcal{S}$ over $\mathbb{M}$ is the
vector bundle such that for every $x \in \mathbb{M}$,
$\mathcal{S}_x$ is the spinor module over the cotangent space
$\mathbf{T}^*_x \mathbb{M}$. At each point $x$, there is therefore
a natural action of $\mathbf{Cl} (\mathbf{T}^*_x \mathbb{M})$ on
$\mathcal{S}_x$, with $\mathbb{C} \otimes \mathbf{Cl}
(\mathbf{T}^*_x \mathbb{M})\simeq \mathbf{End} (\mathcal{S}_x)$;
this action will be denoted by $\mathbf{c}$.

\

On $\mathcal{S}$, there is a canonical elliptic first-order
differential operator called the Dirac operator and denoted
$\mathbf{D}$. In a local orthonormal frame $e_i$, with dual frame
$e_i^*$, we have
\[
\mathbf{D}=\sum_{i} c(e^*_i) \nabla_{e_i},
\]
where $\nabla$ is the Levi-Civita connection. We have an analogue
of Weitzenb\"ock formula which is the celebrated Lichnerowicz
formula (see Theorem 3.52 in \cite{Ber-Ge-Ve}):
\[
\mathbf{D}^2=\Delta +\frac{s}{4},
\]
where $s$ is the scalar curvature of $\mathbb{M}$ and $\Delta$ is
given in a local orthonormal frame $e_i$ by
\[
\Delta=-\sum_{i=1}^d (\nabla_{e_i}\nabla_{e_i}- \nabla_{
\nabla_{e_i} e_i }).
\]

After these preliminaries, we can now turn to the proof of the
local index theorem for $\mathbf{D}$.

\

The first crucial step in all  heat equations approaches of index
theorems is the McKean-Singer type formula (see \cite{McK-Sin} and
Theorem 3.50 in \cite{Ber-Ge-Ve}):
\[
\mathbf{ind} (\mathbf{D})=\mathbf{Str}
(\mathbf{P}_t)=\int_{\mathbb{M}}\mathbf{Str}\text{ } p_t (x,x) dx
, \quad t >0,
\]
where $\mathbf{P}_t=e^{-\frac{1}{2} t \mathbf{D}^2}$, $p_t$ is the
corresponding Schwartz kernel, and $dx$ is the Riemannian volume
form.

By using the results of Section 3 we now show that we actually
have
\[
\lim_{t \rightarrow 0} \mathbf{Str}\text{ } p_t (x,x)dx=\left(
\frac{1}{2i\pi} \right)^{\frac{d}{2}}[A(\mathbb{M})]_d (x).
\]
This last statement is first due to Patodi \cite{Pat} and Gilkey
\cite{Gi} and implies the index theorem. Indeed, $\mathbf{D}^2$ is
a generalized Laplacian, so that the development (\ref{devpt})
holds for $p_t$. With $K=\frac{d}{2}+1$ and $y=x$, it gives
\[
\mathbf{Str}\text{ } p_t (x,x)=\frac{1}{(2\pi t)^{d/2}}
\sum_{k=0}^{d/2} t^k \text{ } \mathbf{Str}\text{ } \Psi_k (x,x)
+O(t),
\]
uniformly in $x \in \mathbb{M}$. The existence of the above limit
at every point $x$ forces $\mathbf{Str}\text{ } \Psi_k (x,x)=0$
for $k<\frac{d}{2}$. The convergence is therefore uniform on
$\mathbb{M}$, and we can pass to the limit in the McKean-Singer
formula.

\

Let us  fix $x_0 \in \mathbb{M}$ once time for all in what
follows. Let $e_i$ be an oriented synchronous local orthonormal
frame centered at $x_0$ with dual frame $e_i^*$. With a cut-off
function, we extend smoothly the vector fields $e_i$ to be zero
outside a neighborhood of $x_0$. At the center $x_0$ of the frame,
we have:

\[
\Delta=-\sum_{i=1}^d \nabla_{e_i} \nabla_{e_i}
\]
and
\[
[\nabla_{e_i},\nabla_{e_j}]=\mathbf{c} \left( D R (e_i,e_j)
\right),
\]
where $R$ is the Riemannian curvature and where $D R (e_i,e_j) \in
\mathbf{Cl}^2 (\mathbf{T}^*_{x_0} \mathbb{M})$ is associated as
above with the skew-symmetric map $R(e_i,e_j)$. These two formulas
only
hold at the point $x_0$, whereas the results of Section 3 require
an expression of $\mathbf{D}^2$ in a neighborhood of $x_0$. On the
neighborhood $\mathcal{O}$ of $x_0$ where $(e_1,...,e_d)$ is an
orthonormal frame, the Lichnerowicz formula reads
\[
-\frac{1}{2} \mathbf{D}^2=\nabla_0 +\frac{1}{2} \sum_{i=1}^d
\nabla_i^2,
\]
where
\[
\nabla_i =\nabla_{e_i}, \text{ }1 \le i \le d, \quad
\nabla_0=\nabla_{e_0} -\frac{s}{8}, \quad e_0=-\frac{1}{2}
\sum_{i=1}^d \nabla_{e_i} e_i.
\]
The vector field $e_0$ vanishes at $x_0$, but in general it does
not vanish in a neighborhood of $x_0$ and it can not be neglected
(for instance, on the two-dimensional unit sphere, a computation
in polar coordinates gives $e_0=-\frac{1}{2} \tan \left(
\frac{r}{2} \right) \frac{\partial}{\partial r}$, where $r$ is the
distance to $x_0$). We are therefore in the framework of Section 3
with
\[
\mathcal{L}=-\frac{1}{2} \mathbf{D}^2, \quad V_i=e_i, \text{ }0
\le i \le d, \quad \mathcal{F}_0=-\frac{s}{8}, \quad
\mathcal{F}_i=0, \text{ }1 \le i \le d.
\]

\ We now apply Corollary \ref{corollary-diag} with $N=d+1$. For
$t >0$ sufficiently small, let $\Theta_t (x_0) \in \mathbf{End}
(\mathcal{S}_{x_0}) \simeq \mathbb{C} \otimes \mathbf{Cl}
(\mathbf{T}^*_{x_0} \mathbb{M})$ be the Clifford element
\[
\Theta_t (x_0) =\mathbb{E} \left( \chi^{d+1}_t (B^{t,b^*_t})
\frac{e^{-\frac{\mid b^*_t \mid^2}{2t}}}{\mid \det \mathbf{J}_t
\mid} \exp \left( \sum_{I, d(I) \le d+1 } \Lambda_I
(B^{t,b^*_t})_t \mathcal{F}_{I}(x_0) \right) \right) ,
\]
where, as in Corollary \ref{corollary-diag}, $b^*_t$ is the
solution of the equation
\[
\sum_{I, d(I) \le d+1 } \Lambda_I (B^{t,b})_t e_{I}(x_0)=0,
\]
$\mathbf{J}_t$ is the differential at $b^*_t$ of the map $b
\rightarrow \exp \left( \sum_{I, d(I) \le d+1 } \Lambda_{I}
(B^{t,b})_t e_{I} \right)(x_0)$, and $\mathcal{F}_I
=\nabla_I-\nabla_{e_I}$, with the same convention as in this
corollary when the cutoff identically vanishes or when there is no
such $b^*_t$. Up to the cutoff and to the Jacobian,
$\Theta_t (x_0)$ has to be thought of as the expectation of the
holonomy $\exp \left( \sum \Lambda_I (B)_t \nabla_{I} \right)$,
conditionally to the fact that the flow $\exp \left( \sum
\Lambda_{I} (B)_t e_{I} \right)$ comes back to $x_0$. From
Corollary \ref{corollary-diag}, we have
\begin{align}\label{dirac-diag}
p_t (x_0,x_0)=\frac{1}{(2\pi t)^{d/2}} \Theta_t (x_0) +O\left(
t^{\frac{1}{2}} \right), \quad t \rightarrow 0.
\end{align}
\begin{proposition}\label{spin-holonomy}
\[
\lim_{t \rightarrow 0} \frac{\mathbf{Str} \text{ } \Theta_t
(x_0)}{t^{d/2}}=\frac{1}{2^{d/2} (d/2)!}\mathbf{Str} \text{ }
\mathbb{E} \left(
 \left. \left(\sum_{1 \le i<j \le
d } D R (e_i,e_j) \int_0^1 B^i_s dB^j_s-B^j_s dB^i_s
\right)^{\frac{d}{2}} \right| B_1=0 \right)
\]
\end{proposition}

\begin{proof}

First, let us observe that due to the scaling property of Brownian
motion, as in the second step of the proof of Theorem
\ref{main-estimation},
\[
\Theta_t (x_0) =\mathbb{E} \left( \bar{\chi}_t
\frac{e^{-\frac{\mid \bar{b}_t \mid^2}{2}}}{\mid \det
\bar{\mathbf{J}}_t \mid} \exp \left( \sum_{I, d(I) \le d+1 }
t^{d(I)/2} \Lambda_I (B^{1,\bar{b}_t})_1 \mathcal{F}_{I}(x_0)
\right) \right) ,
\]
where $\bar{b}_t$ is the solution $b$, $\mid b \mid <3\rho /\sqrt{t}$,
of the equation
\[
\sum_{I, d(I) \le d+1 } t^{\frac{d(I)-1}{2}} \Lambda_I (B^{1,b})_1
e_{I}(x_0)=0,
\]
and where $\bar{\mathbf{J}}_t$ and $\bar{\chi}_t$ are the
corresponding Jacobian and cutoff. The integrand is set to $0$
when the cutoff vanishes identically or when there is no such
$\bar{b}_t$. When the integrand is nonzero, with the notations of
the second step of the proof of Theorem \ref{main-estimation},
$\bar{b}_t=F_{\sqrt{t}}^{-1} (0)$,
$\bar{\mathbf{J}}_t=D_b F_{\sqrt{t}} (\bar{b}_t)$ and
$\bar{\chi}_t=\bar{\chi}_{\sqrt{t}} (\bar{b}_t)$, so that, from
(\ref{jacobian}) and (\ref{inverse-bound}),
\[
\mid \bar{b}_t \mid \le C \sqrt{t} (1+\Xi), \quad \parallel
\bar{\mathbf{J}}_t -\mathbf{Id} \parallel \le \frac{1}{2}.
\]
Moreover, for a fixed loop $B^{1,0}$, since $F_0 (0)=0$ and $D_b
F_0 (0)=\mathbf{Id}$, the implicit function theorem shows that
$\bar{b}_t$ exists as soon as $t$ is small enough. Therefore,
almost surely, when $t \rightarrow 0$,
\[
\bar{b}_t \rightarrow 0, \quad \bar{\mathbf{J}}_t \rightarrow
\mathbf{Id}, \quad \bar{\chi}_t \rightarrow 1.
\]

Let us now rewrite more explicitly the exponential that appears in
$\Theta_t (x_0)$ as a Clifford element.

If $1 \le i , j \le d$, we have in a neighborhood of $x_0$
\[
[\nabla_{i},\nabla_{j}]=\nabla_{[e_i,e_j]}+\mathbf{c} \left( D R
(e_i,e_j) \right)
\]
with
\[
D R (e_i,e_j)=\frac{1}{2}\sum_{1 \le k<l \le d} \langle R(e_i,e_j)
e_k, e_l \rangle e^*_k e^*_l \in \mathbf{Cl}^2 (\mathbf{T}^*
\mathbb{M}),
\]
and $[e_i,e_j](x_0)=0$. Since the Levi-Civita connection is a
Clifford connection, if $1 \le i , j , k \le d$, we have then
\begin{align*}
[\nabla_i,[\nabla_j,\nabla_k]]&=[\nabla_i,\nabla_{[e_j,e_k]} +
\mathbf{c} \left( D R
(e_j,e_k) \right)] \\
 & =\nabla_{[e_i,[e_j,e_k]]} + \mathbf{c} \left(DR (e_i,[e_j,e_k])+ \nabla_i D R
(e_j,e_k) \right).
\end{align*}
and we observe that $DR (e_i,[e_j,e_k])+ \nabla_i D R (e_j,e_k)$
is still a section of $\mathbf{Cl}^2 (\mathbf{T}^* \mathbb{M})$.
In the same way, since $\nabla_0=\nabla_{e_0} -\frac{s}{8}$,
\[
[\nabla_i,\nabla_0]=\nabla_{[e_i,e_0]} + \mathbf{c} \left( DR
(e_i,e_0) \right) -\frac{1}{8} e_i (s).
\]

 More generally, a recurrence procedure shows that for every word
$I$, in a neighborhood of $x_0$,
\[
\nabla_I=\nabla_{e_I}+f_I +\mathbf{c}(\omega_I),
\]
where $f_I$ is a smooth function and $\omega_I$ a smooth section
of $\mathbf{Cl}^2 (\mathbf{T}^* \mathbb{M})$. Therefore
\[
\mathcal{F}_I (x_0)=f_I (x_0) +\mathbf{c}(\omega_I (x_0)),
\]
and for the words $I$ such that $d(I) \le 2$, we have
\[
\omega_{(i)}=0, \text{ } 0 \le i \le d, \quad \omega_{(i,j)}
(x_0)=D R (e_i,e_j), \text{ } 1 \le i,j \le d.
\]
Since scalars commute with the Clifford elements, we deduce that
\[
\exp \left( \sum_{I, d(I) \le d+1 } t^{d(I)/2} \Lambda_I
(B^{1,\bar{b}_t})_1 \mathcal{F}_{I}(x_0) \right)=e^{X_t} \text{ }
\mathbf{c} \left(\exp \left( \Omega_t \right) \right) ,
\]
where
\[
X_t=\sum_{I, d(I) \le d+1 } t^{d(I)/2} \Lambda_I
(B^{1,\bar{b}_t})_1 f_I (x_0)
\]
is a scalar term and
\[
\Omega_t=\sum_{I, 2 \le d(I) \le d+1 } t^{d(I)/2} \Lambda_I
(B^{1,\bar{b}_t})_1 \omega_I (x_0) \in \mathbf{Cl}^2
(\mathbf{T}^*_{x_0} \mathbb{M}).
\]
Therefore
\[
\mathbf{Str} \text{ }\Theta_t (x_0)= \mathbb{E} \left(
\bar{\chi}_t \frac{e^{-\frac{\mid \bar{b}_t \mid^2}{2}}}{\mid \det
\bar{\mathbf{J}}_t \mid} e^{X_t} \text{ } \mathbf{Str} \text{
}\exp \left( \Omega_t \right) \right).
\]
But now, since $\Omega_t \in \mathbf{Cl}^2 (\mathbf{T}^*_{x_0}
\mathbb{M})$, according to (\ref{supertraceformula-clifford}), for
any $k<\frac{d}{2}$,
\[
\mathbf{Str} \text{ } \Omega_t^k =0.
\]
Moreover, the elements of $\mathbf{Cl}^2 (\mathbf{T}^*_{x_0}
\mathbb{M})$ act on the spinor module by skew-Hermitian
endomorphisms, so that $\exp \left( \tau \Omega_t \right)$ acts by
an isometry. From Taylor formula, we get therefore
\[
\left| \mathbf{Str} \text{ }\exp \left( \Omega_t \right)
-\frac{1}{(d/2)!} \mathbf{Str} \text{ } \Omega_t^{\frac{d}{2}}
\right| \le C \parallel \Omega_t \parallel^{\frac{d}{2}+1}.
\]
On the support of the localized integrand, $X_t$ is bounded by
(\ref{cutoff2}), and from (\ref{polynomial}), $\parallel \frac{1}{t}
\Omega_t \parallel \le C (1+\Xi)^L$. Finally, almost surely, when
$t \rightarrow 0$, $X_t \rightarrow 0$ and
\[
\frac{1}{t} \Omega_t \rightarrow \sum_{1 \le i,j \le d}
\Lambda_{(i,j)} (B^{1,0})_1 D R (e_i,e_j)=\frac{1}{2}\sum_{1 \le
i<j \le d } D R (e_i,e_j) \int_0^1 B^{1,0,i}_s dB^{1,0,j}_s
-B^{1,0,j}_s dB^{1,0,i}_s .
\]
Since $\Xi$ has moments of all orders, we conclude by dominated
convergence that
\[
 \lim_{t \rightarrow 0} \frac{\mathbf{Str} \text{ } \Theta_t
(x_0)}{t^{d/2}}=\mathbb{E} \left( \mathbf{Str} \text{ }
\frac{1}{(d/2)!} \left(\frac{1}{2}\sum_{1 \le i<j \le d } D R
(e_i,e_j) \int_0^1 B^{1,0,i}_s dB^{1,0,j}_s -B^{1,0,j}_s
dB^{1,0,i}_s \right)^{\frac{d}{2}} \right).
\]
This is the expected result because the law of the loop $B^{1,0}$
is the law of $(B_s)_{0 \le s \le 1}$ conditioned by $B_1=0$.
\end{proof}

We can now obtain the required limit:

\begin{theorem}
\[
\lim_{t \rightarrow 0}\mathbf{Str} \text{ }p_t
(x_0,x_0)dx_0=\left( \frac{1}{2i\pi}
\right)^{\frac{d}{2}}[A(\mathbb{M})]_d (x_0).
\]
\end{theorem}

\begin{proof}
By the estimate (\ref{dirac-diag}), which comes from Theorem
\ref{main-estimation} and Corollary \ref{corollary-diag}, and by
Proposition \ref{spin-holonomy}, we get
\begin{align*}
& \lim_{t \rightarrow 0}\mathbf{Str} \text{ }p_t (x_0,x_0)dx_0 \\
=& \frac{1}{(4\pi)^{d/2}(d/2)!}\mathbf{Str} \text{ } \mathbb{E}
\left(
 \left. \left(\sum_{1 \le i<j \le
d } D R (e_i,e_j) \int_0^1 B^i_s dB^j_s-B^j_s dB^i_s
\right)^{\frac{d}{2}} \right| B_1=0 \right) e^*_1 \wedge ...
\wedge e^*_d .
\end{align*}
Let us denote by $\Omega_{ij}=\langle \Omega e_i , e_j \rangle$ the
entries of the curvature form $\Omega$, which are two-forms. Due
to the symmetries of the Riemannian curvature,
\[
\Omega_{ij}=\sum_{1 \le k<l \le d} \langle R(e_i,e_j) e_k, e_l
\rangle e^*_k \wedge e^*_l.
\]
From (\ref{supertraceformula-clifford}) and from the expression of
$D R (e_i,e_j)$, the above expression is therefore also equal to
the $d$-form piece of
\[
\frac{1}{(2i\pi)^{d/2}(d/2)!} \mathbb{E} \left(
 \left. \left(\sum_{1 \le i<j \le
d } \frac{1}{2} \Omega_{ij} \int_0^1 B^i_s dB^j_s-B^j_s
dB^i_s \right)^{\frac{d}{2}} \right| B_1=0 \right).
\]
This last expression is also the $d$-form piece of
\[
\frac{1}{(2i\pi)^{d/2}} \mathbb{E} \left(
 \left. \exp \left(\sum_{1 \le i<j \le
d } \frac{1}{2} \Omega_{ij} \int_0^1 B^i_s dB^j_s-B^j_s
dB^i_s \right) \right| B_1=0 \right),
\]
which turns out to be the $d$-form piece of the $A$-genus
\[
\frac{1}{(2i\pi)^{d/2}}A(\mathbb{M})=\frac{1}{(2i\pi)^{d/2}}\det
\left( \frac{ \Omega}{2 \sinh \frac{1}{2} \Omega}
\right)^{\frac{1}{2}},
\]
from L\'evy's area formula (see for instance Lemma 7.6.6 in
\cite{Hsu}).
\end{proof}

Let us finally come back to the equality of limits that was
announced in the introduction. From (\ref{kernel-formula}) and
from the proof of Corollary \ref{corollary-diag}, for every $N \ge
1$, $(2\pi t)^{d/2} p^N_t (x_0,x_0)$ is the expectation that
appears in this corollary. For $N=d+1$, this is $\Theta_t (x_0)$.
For $N=2$, since $e_0 (x_0)=0$ and $[e_i,e_j](x_0)=0$, the
equation that defines $b^*_t$ simply reads $b=0$, and the proof of
Proposition \ref{spin-holonomy} applies in the same way. Therefore
\[
\lim_{t \rightarrow 0} \mathbf{Str} \text{ } p^{2}_t
(x_0,x_0)=\lim_{t \rightarrow 0} \mathbf{Str} \text{ } p^{d+1}_t
(x_0,x_0)=\lim_{t \rightarrow 0} \mathbf{Str} \text{ } p_t
(x_0,x_0).
\]

\end{document}